\documentclass[10pt,letterpaper]{article}
\usepackage{t1enc}
\usepackage[latin1]{inputenc}
\usepackage[english]{babel}
\usepackage{graphics}
\begin{document}

\title{Isochronous centers of Lienard type equations and Applications}
\author{A. Raouf Chouikha 
\footnote
{Universite Paris 13 LAGA, Villetaneuse 93430. {\small chouikha@math.univ-paris13.fr}  }
}
\date{}
\maketitle
{\it Dedicated to the 60th anniversary of Professor Jean-Marie Strelcyn}

\bigskip

\begin{abstract}
In this work we study the equation \  $(E)\ \ddot x + f(x) \dot x^2 + g(x) = 0$ \ with a center at $0$ and investigate conditions of its isochronicity. When $f$ and $g$ are analytic (not necessary odd) a necessary and sufficient condition	for the isochronicity of $0$ is given. This approach allows us to present an algorithm for obtained conditions for a point of (E) to be an isochronous center. 
In particular, we find again by another way the isochrones of the quadratic Loud systems \ $(L_{D,F})$. Some classes of Kukles are also considered. Moreover, we classify a 5-parameters family of reversible cubic systems with isochronous centers.

 {\it Key Words and phrases:} \ period function, monotonicity, isochronicity, center, polynomial systems\footnote
{2000 Mathematics Subject Classification  \ 34C25, 34C35}. 
\end{abstract}

\bigskip

\section{ Introduction and statement of results}

Consider the planar differential system $$ \dot x= P(x,y),\qquad  \dot y = Q(x,y),\eqno(1)$$ where \ $\dot x = \frac {dx}{dt}, \dot y = \frac {dy}{dt}, 
 \ P(x,y)$\ and \ $Q(x,y)$\ are analytic functions defined in an open subset of \ $R^2$.\\
 To study the integrability of system (1) we may investigate local first integrals.\\
 Consider the case where \ $P$\ and \ $Q$\ can be written 
 $$P(x,y) = -y + U(x,y) , \qquad Q(x,y) = x + V(x,y) \eqno(2)$$ where \ $U$\ and \ $V$ \ are convergent series without linear terms.\\
 Recall that system (2) has a center at the origin $0$ if and only if it has a Lyapounov first integral. Moreover, the center is isochronous if and only if system (2) is linearizable. Namely, the origin $0$ of system (2) is a center if all orbits in a neighborhood are closed. $0$ is an isochronous center if the period of oscillations is the same for all these orbits.\\
 We say that system (2) is linearizable if there is analytic change of coordinates in the neighborhood of $0$ bringing the system into the linear one.\\
 Using the complex variable \ $u = x + iy$\ we may write (2) as a equation 
 $$\dot u = iu + R(u,\bar{u}) \eqno(3)$$ where \ $R(u,\bar{u})$\ is an analytic function.\\
 To better understand isochronicity phenomena we substitute 
 $$r^2 = u \bar{u} , \qquad \theta = \arctan (\frac {Im\ u}{Re\ u}).$$
 we get $$\frac {dr}{d\theta } = ir \frac {\dot u\bar{u}+u\bar{u}}{\dot u\bar{u}-u\bar{u}} = \frac {2r(R \bar{u}+\bar{R}u)}{2r-R \bar{u}+\bar{R}u} \eqno(4)$$
 Let \ $r(\theta,\rho)$\ denotes the solution of (4) verifying \ $r(0,\rho) = \rho.$\\
 The quantity \ $ r(2\pi,\rho)$\ is the return map starting from the $\theta = 0$ axis.\\
 Write $$r(\theta,\rho) = \rho + u_2(\theta)\rho ^2 + u_3(\theta)\rho ^3 +.... \quad with \quad u_k(0) = 0.$$
 The origin is a center of equation (3) if \ $r(2\pi,\rho) = \rho$\ namely \ $u_k(2\pi) = 0.$\\
 When system (3) has a center define the period function \ $T(\rho)$\ as the time spent by the closed orbit to return around the origin.\\
Let us now write the angular speed
 $$\frac {d\theta}{dt} = - i\frac {\dot x\bar{x}+x\bar{\dot x}}{2 r^2}.$$
 So, near the origin the period function may be expressed   $$T(\rho) = 2\pi + \sum_{1\leq k}t_k(2\pi)\rho^k.$$
 Then, the center is isochronous if and only if  \ $t_k(2\pi) = 0$\ for all \ $k \geq 1.$\\
 It is known that the first $k$ verifying \ $t_k(2\pi) \neq 0$\ is necessary an even number.\\
Notice that this number called {\it period quantity} plays a role in problems of bifurcations of critical periods.\\

Isochronicity phenomena has been widely studied not only for its impact in stability theory, but also for its relationship with bifurcation and boundary value problems.\\
Up to now the center problem as well as the isochronicity problem for systems of the form (2) is solved in many cases of the systems with homogeneous quadratic and cubic nonlinearities.\\
 
For the center $0$  of (2) the largest neighborhood of $0$  which is covered by periodic orbits is  the {\it period annulus} of $0$ denoted by $\gamma _0$. A center is said  to be a {\it global center} when its period annulus is the whole plane.\\ Let the period function $T$ associated to any periodic orbit $\gamma $ in $\gamma _0$. A center is said to be {\it nondegenetate} when the linearised vector field at the critical point has two nonzero eigenvalues. It is well known that only nondegenerate centers can be isochronous.\\
 When the differential system is analytic (that means \ $P(x,y)$\ and \ $Q(x,y)$\ are analytic functions) it implies that the period annulus of an isochronous center is unbounded. \\

 The orbits may be parameterized for example by choosing their initial values in the segment $(0, \pi)$ on the $x$ - axis.\\  Let  $T  : \gamma _0 \rightarrow R $ ,\  be the function defined  by associating to every point  $(x,0) \in \gamma _0$ the minimum period of the trajectory starting at $(x,0)$, to reach the negative x-axis. \ $T$ is the period function and is constant on cycles. We say that $T$ is (strictly) increasing if, for every couple of cycles $\gamma _1$ and $\gamma _2 ,\  \gamma_1$ included in $ \gamma _2$ , we have \ $T(\gamma _1) \leq T(\gamma _2)$ \  ( $T(\gamma _1) < T(\gamma _2) $).\\ 
$0$ is an isochronous center if $T$ is constant in a neighborhood of $0$.\\

One of the most studied systems are those of the form 
$$ \dot x = y,\qquad  \dot y = - Q(x,y),\eqno(5)$$
which are equivalent to the second order equation
$$\ddot x + Q(x,\dot x) = 0.$$
Notice that a stationary point of a second order differential equation is a center if it corresponds to a center for the equivalent planar system.\\
Many second order differential equations arising from mechanics and electricity can be reduced to that of suitable systems. Such systems are also called Kukles systems when \ $Q(x,y)$\ is a polynomial with real coefficients of degree $d$ without $y$ as a divisor.\\

  Volokitin and Ivanov proved that for every positive integer $m$, the origin $0$ is an isochronous center of the equation 
$$\ddot x + x^{2m-1} \dot x (x^2+\dot x^2) + x = 0.$$
More precisely, they proved that the equivalent system 
$$ \dot x = y,\qquad  \dot y = - x - x^{2m-1} y (x^2+y^2) + x,$$
has only trivial polynomial commutators.
\\

The Lienard equation 
$$\ddot x + f(x) \dot x + g(x) = 0$$
and the isochronicity of its center has motivated many authors. \\ In particular, Christopher, Devlin, Lloyd and Sabatini [C-S] proved if $f$ and $g$ are analytic odd
functions of $x$ with $g'(0) = 1$ and $x g(x)>0$ in a neighborhood of the origin. Then 
the Lienard equation has an isochronous center at\ $0$\ if and only if 
$$ 
g(x) = x + {1\over x^3}\biggl(\int_0^x \xi f(\xi)\,d\xi\biggr)^2.
$$
\bigskip
 
Here is another interesting equation of Lienard type
$$(E) \qquad \ddot x + f(x) \dot x'^2 + g(x) = 0$$
or its equivalent system denoted also $(E)$ $$\dot x = y , \qquad \dot y = -g(x) - f(x) y^2$$
where $f, g$ are of class $C^1$ in a neighborhood of $0$ verifying the condition $x g(x) > 0$ for $x \neq 0$ and $\frac {dg}{dx}(0) > 0$.
A such equation has special properties.\\ In particular, it can be reduced to a conservative equation and then has a first integral. Opposite to the Lienard equation which cannot be reduced (in general) to a conservative equation  and a first integral is unknown (except for very few cases).\\

Equation (E) has many physical applications, [L-R]. Indeed, that is a model of one dimensional oscillator studied at the classical and also at the quantum level. The following particular case interested specially the physicists
$$f(x) = \frac {- \lambda x}{1+\lambda x^2}, \qquad g(x) = \frac {\alpha ^2x}{1+\lambda x^2}$$
The general solution takes the form \ $ x(t) = A \sin (\frac {2\pi}{T(A)}t + \phi).$
Curiously, that is the only case which permits to explicitly determine the amplitude dependance of the period function 
$$T = T(A) = \frac {2\pi \sqrt{1 + \lambda A^2}}{\alpha}.$$ 
$T(0) = \frac {2\pi}{\alpha}$ \quad corresponds to the center $0$. We see that $T(A)$ is an increasing function.  \footnote{see below in Appendice 1 some facts concerning monotonicity properties which may interest physicists.}
\newpage
 When $f$ and $g$ are of class $C^1$ in $(E)$, Sabatini [S] gave a sufficient condition for the monotonicity of the period $T$ or for the isochronicity of $0$. Such a condition is also necessary when $f$ and $g$ are odd and analytic.\\

This paper is organized as follows. 
We give here a sufficient condition for the monotonicity of the period $T$ of equation $(E)$. When $f$ and $g$ are  analytic - not necessary odd - we establish a necessary and sufficient condition for the isochronicity of the center $0$. We then extend some Sabatini's results.\\ This fact allows us to present an algorithm for finding conditions for a critical point to be an isochronous center of $(E)$. It is based on a transformation of equation $(E)$ to a conservative one and in using an Urabe theorem. \\ Applying to quadratic systems, we
find another isochronicity condition for the dehomogenized Loud systems \ $(L_{D,F})$ :
$$(C_2)\quad 4F^3 + 24DF + 24D^2 + 2DF^2 - F^2 - 4F - 2D + 1 = 0.$$
Combining with the classical relation 
$$(C_1)\quad 4 F^2 + 10 DF + 10 D^2 - D - 5F +1 = 0,$$
it yields exactly the four isochrones for \ $(L_{D,F})$: say $(L_{0,1}), (L_{\frac {-1}{2},2}), (L_{0,\frac {1}{4}}), (L_{\frac {-1}{2},\frac {1}{2}}).$\\  
  Another application concerns the monotonicity of the period function for the reduced Kukles systems.\\ Finally, we contribute to the study of some (non homogenous) cubic systems with an isochronous center. We prove that the 5-parameters system
 $$ (C) \qquad \left\{ \begin{array}{c} \dot x = - y + b x^2y\\
 \dot y = x + a_1x^2 + a_3y^2 + a_4x^3 + a_6xy^2 \end{array}\right. $$
 admits only four classes of isochrone systems.
 This produces examples of cubic reversible systems verifying $R_3 = \frac {a1-a3}{4}\neq 0$ that are not covered by the classification of Chavarriga and Garcia ([C-S] sect. 12).\\
 
 \bigskip
 
 {\it Aknowledgment }: I would like to thank Javier Chavarriga whose his help permits us to state below Theorem 4-3 . 
 
\newpage

\section {Isochronous centers of the system (E)}
Let us consider $$(E)\quad \dot x = y , \qquad \dot y = -g(x) - f(x) y^2$$ and the integrals
$$
F(x) = \int_0^x f(s) ds , \qquad \phi (x) = \int_0^x e^{F(s)} ds.\qquad (6)
$$
   If $f,g$ are of class $C^1$ in a neighborhood $N_0$ of $0$, then the function $u = \phi (x)$ is invertible in $N_0$. The following in particular allows us to transform  $(E)$ into a conservative one, ( Lemma 1 of [S])
   $$\ddot u + \tilde g(u) = 0$$

\bigskip

{\bf Lemma 2-1}\quad {\it Under the above hypothesis $x(t)$ is a solution to  
$$(E) \qquad \ddot x + f(x) \dot x^2 + g(x) = 0$$
if and only if \ $u(t) = \phi (x(t))$\ is a solution to  
$$(E_c)\qquad \ddot u + g(\phi ^{-1}(u)) e^{{F(\phi ^{-1}(u))}} = 0$$
and the function \ $\tilde g(u) = g(\phi ^{-1}(u)) e^{{F(\phi ^{-1}(u))}}$\ is of class \ $C^1$.\\ Moreover, if $0$ is a center of $(E)$ then it is also a center of $(E_c)$}

\bigskip

Consider the change \ $u = \phi (x).$\ It gives \ $\dot u = e^{F(x)}\dot x$\ and $$ \ddot u= f(x)e^{F(x)}\dot x^2 + e^{F(x)}\ddot x = e^{F(x)}[\ddot x + f(x) \dot x^2].$$ Then 
$$\ddot u= - g(x)e^{F(x)} = - g(\phi ^{-1}(u))e^{F(\phi ^{-1}(u))}.$$
Moreover, $(E)$ has a center if \ $x g(x) > 0.$ in $N_0$.\ Or equivalently  \ $\phi (x) g(x) > 0 $\ since by definition \ $x \phi (x) > 0$ in $N_0$.\
 It implies \ $\phi (x) g(x) e^{F(x)} > 0$.\\ That means $(E_c)$ has a center at $0$ since \ $u \tilde g(u) > 0$ for $u$ closed to $0$ and $u \neq 0$.\\
 Notice that when \ $g$ \ is analytic \ $x g(x) > 0$\ in \ $N_0$\ is a necessary and sufficient condition for the origin $0$ to be a center.\\   

By this lemma we may deduce at first, trivial cases of isochronicity for equation (E). Taking  
$$
\tilde g(u) = g(\phi ^{-1}(u)) e^{F(\phi ^{-1}(u))} = K u
$$
 with $K > 0$.\\ We see that when\  $F \in C^1(N_0)$ \ all equations of the form 
$$
\ddot x + F'(x) \dot x^2 + e^{-F(x)} \int_0^x e^{F(s)} ds = 0\qquad (7)
$$
have a isochronous center at $0$ ( the symbol $'$ means $\frac {d}{dx}.$)\\
There are many others equations reducing to the trivial linear one, as we will see below .\\ 

 Consider the conservative system 
 $$(E_c) \qquad \dot u = v , \qquad \dot v = - \tilde g(u)$$ 
 where \ $\tilde g$\ is \ $C^1$\ and such that \ $g(0) = 0, \ \frac {dg}{du}(0) = 1.$ \ Let the integral  
 $$\tilde G(u) = \int_0^u \tilde g(s) ds $$
 The following proved by Urabe permits to characterize isochronous centers for systems $(E_c)$, see [U].

 \bigskip
 
 {\bf Lemma 2-2}\quad {\it Let \ $\tilde g(u) $\ be a \ $C^1$ \  function defined in $N_0$ a neighborhood of $0$ verifying $u \tilde g(u) > 0$ in $N_0$. Then the system $(E_c)$ has an isochronous center at the origin $0$ if and only if by the transformation  \ $X^2 = 2\tilde G(u)$\ where \ $\frac {X}{u} > 0,\ \tilde g(u)$\ may be written 
 $$\tilde g(u) = \frac {X}{1 + h(X)}$$
 where \ $h(X)$\ is a \ $C^1$ \ odd function}

 \bigskip
 When $h(X)$ is non trivial we will call it in the sequel the Urabe function.\\
A simple example of isochronous case may be obtained in solving the differential equation
$$
 S_g = 5 \tilde g''^2(u) - 3 \tilde g'(u) \tilde g'''(u) = 0.\qquad (8)
$$
 - It is known that  \quad $S(\tilde g) = 5 \tilde g''^2(u) - 3 \tilde g'(u) \tilde g'''(u) \neq 0$\quad  implies the monotonicity of the period function for System $(E_c)$ near a center. This criteria has been introduced by R. Schaaf, [Sc].\\
 Resolution of $(3)$ gives 
 $$\tilde g(u) =  1 - (1 + 2u)^{-\frac{1}{2}}$$ 
 which corresponds to the odd trivial function \ $h(X) = X$\ and the corresponding isochronous potential is 
 $$\tilde G(u) = 1 + u - \sqrt {1 + 2u}$$
 where \ $-\frac{1}{2} < u < \frac{1}{2}$\ so that this potential is analytic.\\
 
Using preceding lemmata one may deduce other (non trivial) classes of equations $(E)$ having an isochronous center at $0$.

\newpage

{\bf Proposition 2-3}\quad {\it Let \ $f,g \in C^1(N_0)$,\ where $N_0$ is a neighborhood of $0$ and $x g(x) > 0$\ in \ $N_0$\ then the origin $0$ is a center of Equation 
$$(E)\qquad \ddot x +f(x) \dot x^2 + g(x) = 0$$
 (if $g$ is analytic then $0$ is a center if and only if $x g(x) > 0$ in $N_0$).\\ 
Suppose the integral of $f$ : $F(x) = \int_0^x f(t) dt $ verifies  $$e^{F(x)} + e^{F(-x)} = 2\quad  
and \quad g(x) = K \frac {x}{e^{2F(x)}}$$ where $K$ is a positive constant, then \ $0$\ is an isochronous center of (E).\\
In particular if in addition $g$ is odd, then $0$ is an isochronous center of (E) if and only if \ $f(x) \equiv 0 $\  and \ $ g(x) = K x$.}   

\bigskip
Notice at first that the classes of functions $g$ given by Proposition 2-3 is different as (7). Indeed, suppose
$$  \frac {x}{e^{2F(x)}} = e^{-F(x)} \int_0^x e^{F(s)} ds$$
a quick calculation shows that no function $F \neq 0$ verifies this functional equation ( $F\equiv 0$, corresponds to an odd function $g$).
 
  \bigskip
  
  {\bf Proof}\quad Let $X$ defined by $$\frac {1}{2} X^2 = \int_0^x g(s) e^{2F(s)} ds.$$
  Consider the following change of variable
  $$u = X + H(X)$$
  where \ $H(X)$\ is an even \ $C^2$\ function defined in a neighborhood of $0$ such that $H(0) = 0 $ and $h(X)$ is its derivative such that \ $\frac {X}{u} > 0$.\  Since \    $u = \phi (x)$\ then  \ $\phi (x) = X + H(X)$\ is invertible in a neighborhood of $0$. \\
  Let \ $\tilde f(x) = e^{F(x)} - 1$.\ Then, \ $ \frac {d\tilde f}{dx} =  f (1 + \tilde f)$.\\
  So, 
  $$f = \frac {\frac {d\tilde f}{dx}}{1 + \tilde f} \quad and \quad g = K x e^{-2F} = \frac {K x}{(1 + \tilde f)^2}$$
  Thus,  $(E)$ is equivalent to 
$$
 \quad  \ddot x + \frac {\frac {d\tilde f}{dx}}{1 + \tilde f} \dot x^2 + \frac {K x}{(1 + \tilde f)^2} = 0.\qquad (9)
$$
   However, \ $e^{F} = 1 + \tilde f$\ implies \ $u = \phi (x) = \int_0^x e^{F(s)} ds = x + \tilde F(x)$\ where \ $\tilde F(x) = \int_0^x \tilde f(s) ds$\ and \ $\tilde g(u) = g(x) e^{F(x)} = \frac {K x}{1 + \tilde f(x)}$\\
   Moreover,
   $$\frac {1}{2} X^2 = \int_0^x g(s)e^{2F(s)} ds = \frac {1}{2} K x^2$$
   Then ,$$X = x \sqrt K, \quad and \quad \tilde g(u) = \frac {X\sqrt K}{1 + \tilde f(\frac {X}{\sqrt K})}.$$  
  Therefore, if \ $\tilde f(\frac {X}{\sqrt K})$\ is an odd function or equivalently 
  $$e^{F(x)} + e^{F(-x)} = 2.$$
  Then by Lemma 2-2 (Urabe theorem) Equation $(4)$ has an isochronous center at $0$. When \ $g = g(x) $\ is odd then a necessary and sufficient condition for $0$ to be an isochronous center is : \ $f(x) \equiv 0$\ and \ $g(x) = K x$.\\ 
  
 Consider the following     

 \bigskip
     
  {\bf Lemma 2-4}\quad {\it Equation (E) has an isochronous center at the origin $0$ if and only if 
$$
  (1+ h(X)) \ddot X + \frac {dh}{dX} \dot X^2 + \frac {X}{1+h(X)} = 0\qquad (10)
$$
  where \ $$X^2 = 2 \tilde G(x), \quad \tilde G(x) = \int_0^{x} g(s) e^{2F(s)} ds$$
  has an isochronous center at the origin $0$ where $h(X)$ is an odd function and $h \in C^1(N_0, R), \ N_0$ is choosen such that $1 + h(X) > 0$.}
 
 \bigskip
Using the change\ $X = x \sqrt K$ equation $(10)$ is obviously derived from $(9)$ which is equivalent to $(E)$.\ This lemma implies that $(E)$ has also an isochronous center at the origin.\\ This completes the proof of Proposition 2-3.\\
 
 {\bf Proof of Lemma 2-4}\\  

Condition $x g(x) > 0$ in a neighborhood of $0$ implies the origin is a center of (E).
For $x$ closed to $0$ one has \ $X\frac {X}{(1 + h(X))^2} > 0$. \ This means $0$ is a center of (5).
 Let $H(X) = \int_0^X h(s) ds$ be the primitive of $h(X)$ and denote as above by $u = X + H(X), \ u = \phi (x) = e^{F(x)}.$ \\ Then, $\frac {du}{dt} = (1 + h(X)) \dot X$ and $ \frac {d^2u}{dt^2} = (1 + h(X)) \ddot X +  \frac {dh}{dX} \dot X^2$. Moreover, since  $H(X)$ is an even function and $1 + h(X) \neq 0$ then the relation $u = X + H(X)$ is invertible in a neighborhood of $0$. Let the inverse \ $X = \tilde \phi (u).$\\
 
By Urabe theorem (Lemme 2-2) $(E)$ has an  isochronous center at the origin $0$ if and only if the equation 
 $$\ddot u + \tilde g(u) = 0$$ has   an  isochronous center at the origin $0$, where 
 $$\tilde g(u) = g(x) e^{F(x)} .$$
 Moreover, 
 $$\tilde g(u) = \frac {X}{1+h(X)} = \frac {\tilde \phi(u)}{1+h(\tilde \phi (u))} $$
where \ $h(X)$\ is an odd function. It is also equivalent to assert that the center $0$ of (5) is isochronous.\\

 furthermore, we are able to produce the exact expression of the period function $T$ of (E). This one naturally depends on $h(X)$ defined above.\\ Notice that the Urabe function $h(X)$ also plays a role for the monotonicity of the period function $T$ near the center $0$ of  $(E)$. 
 
 \bigskip
 
  {\bf Proposition 2-5}\quad {\it Let $f,g$ be analytic function, and consider Equation $$(E) \qquad \ddot x + f(x)\dot x^2+g(x) = 0$$ with a center at the origin $0$.\\ Let $X$ defined by $\frac {1}{2} X^2 = \int_0^x g(s) e^{2F(s)} ds$ \ and \ $H(X)$ is such that \\ $ \int_0^x e^{F(s)} ds = X + H(X)$\ and \ $X \int_0^x e^{F(s)} ds > 0$.\\ Then the period function $T$ may be expressed under the form
  $$T = T(c) = 2 \int _{\frac {-\pi}{2}}^{\frac {\pi}{2}} [1 + \frac {dH}{dX}(\sqrt {2c}\sin \theta)] d\theta $$
  where the constant \ $c = \frac {1}{2}\sqrt {X(t=0)}$.\\
  Moreover the following holds :\\ 
(i) \ $\frac {d^3 H}{dX^3}(0) > 0$\ then the period function  $T$ increases in a neighborhood of $0$.\\ 
(ii) \ $\frac {d^3 H}{dX^3}(0) < 0$\ then the period function  $T$ decreases in a neighborhood of $0$.\\
(iii) \ $\frac {d^3 H}{dX^3}(0) = 0$\ is a necessary condition for the center to be isochronous  }  

   \bigskip
   
   {\bf Proof}\quad By Lemma 2-1, $(E)$ is equivalent to the conservative system \\
   \ $\dot u = v, \quad \dot v = - \tilde g(u)$\  where \ $\tilde g(u) = g(x) e^{F(x)}$\ 
   where \ $u = \int_0^x e^{F(s)} ds.$\\
   Then a such function may be written  
   $$\tilde g(u) =  \frac {X}{1+h(X)}$$
   where \ $h(X) = \frac {dH}{dX}.$\\
   
   A calculus yields
   $$\tilde g'(u) = \frac {d\tilde g}{du} = \frac {dX}{du} [\frac {1}{1+h(X)} - \frac {h'(X)X}{(1+h(X))^2}] = \frac {1}{(1+h(X))^2} - \frac {h'(X)X}{(1+h(X))^3}$$
   $$\tilde g''(u) = \frac {-3h'- h''X}{(1+h(X))^4} + \frac {3h'^2X}{(1+h(X))^5}$$
   $$\tilde g^{(3)}(u) = \frac {-4h''- h^{(3)}X}{(1+h(X))^5} + \frac {15h'^2+10h'h''X}{(1+h(X))^6} - \frac {15h'^3X}{(1+h(X))^7}$$
   So,  $$\tilde g'(0) = 1, \quad \tilde g''(0) = - 3h'(0), \quad \tilde g^{(3)}(0) = - 4h''(0) + 15 h'^2(0)$$
   We thus obtain 
   $$5 \tilde g''^2(0) - 3 \tilde g'(0) \tilde g^{(3)}(0) = 12 h''(0) = 12\frac {d^3 H}{dX^3}(0).$$
   By Schaaf criteria [Sc] for the monotonicity of the period function for a conservative system  (see above Equation (8)), 
we easily deduce assertions (i), (ii) and (iii). \\ 

Turning now to the expression of the period function.\\
It is wellknown the period function of the conservative system 
$$(E_c) \qquad \dot u = v , \qquad \dot v = - \tilde g(u)$$ 
with a center at $0$ may be expressed  
$$T(c) = \sqrt {2} \int_{a}^{b} \frac {du}{\sqrt {c - \tilde G(u)}}$$
where \  $\tilde G(u) = \int_0^u \tilde g (s) ds,$\ the constants $a, b$ are such that $a < 0 < b$ and $\tilde G(a) = \tilde G(b) = c$.\\
Recall at first the relation \ $u = X + H(X)$\ is invertible in a neighborhood of $0$. 
Its inverse will be used for change of variables\ $X = \tilde \phi (u)$\ which  
transforms the closed orbits into circles centered at the origin. Then :
$$T(c) = \sqrt {2} \int_{-\sqrt{2c}}^{\sqrt{2c}} \frac {X dX}{\tilde g(u(X))\sqrt {c -\frac {X^2}{2}}}.$$
Finally, another change of variables \ $X = \sqrt{2c} \sin \theta$\ gives 
$$T(c) = 2 \int_{-\frac {\pi}{2}}^{\frac {\pi}{2}}\frac {\sqrt{2c} \sin \theta}{\tilde g(u(\sqrt{2c} \sin \theta))} d\theta $$
$$T(c) = 2 \int _{\frac {-\pi}{2}}^{\frac {\pi}{2}} [1 + \frac {dH}{dX}(\sqrt {2c}\sin \theta)] d\theta .$$

 \bigskip

 As a consequence of Proposition 2-5, we easily deduce that a necessary condition for the center $0$ to be isochronous is \ $h''(0) = 0$. \\ 
 
 In fact, we may obtain a better result as we will see in the sequel : \ condition \ " $h(X)$\ odd "\ is a necessary and sufficient condition for the center $0$ to be isochronous, without supposing $f,g$ odd.\\
 
 Our main result is the following

   \newpage
     
{\bf Theorem 2-6}\quad {\it Let $f,g$ be analytic function in a neighborhood $N_0$ of $0$ , 
and $x g(x) > 0$ for $x \neq 0$ \ 
then  $$(E)\qquad \ddot x + f(x)\dot x^2 + g(x) = 0$$ has an isochronous center at the origin $0$ if and only if 
$$\frac {X}{1+h(X)} = g(x) e^{F(x)}$$ where \ $X$ \ is defined by\ $\frac {1}{2} X^2 = \int_0^x g(s) e^{2F(s)} ds$ \ and \ $h(X)$ \ is an odd function such that \ $ \phi(x) = \int_0^x e^{F(s)} ds = X + \int_0^X h(t) dt$\ and \ $\frac {X}{\phi (x)} > 0$.\\
In particular, when \ $g$ and $f$ are odd then \ $0$\ is an isochronous center if and only if \  $g(x) = e^{-F(x)}\phi (x)$\ (or equivalently \ $h(X) \equiv 0$ ).}

  \bigskip
  
 {\bf Proof} \quad  Theorem 2-6 may be deduced from preceding Lemmata.\\
  Let \ $u = \phi(x) = \int_0^x e^{F(s)} ds$\ and define 
 $$\tilde g(u) = g(x) e^{F(x)}.$$
 Then, by Lemma 2-1 system \ $ (E)$\  is equivalent to the conservative one 
 $$(E_c)\qquad \dot u = v,\qquad \dot v = - \tilde g(u)$$
 where $\tilde g $  is such that \ $\tilde g(0) = 0, \ \frac {d\tilde g}{du}(0) = 1$. Let the integral  $$\tilde G(u) = \int_0^u \tilde g(s) ds.$$

 Moreover, $(E)$ has a center if \ $x g(x) > 0$ in $N_0$\ i.e.  \ $\phi (x) g(x) > 0 $\  since \ $x \phi (x) > 0$ in $N_0$.\
 It implies \ $\phi (x) g(x) e^{F(x)} > 0$.\\ We deduce that $(E_c)$ has a center at $0$ since \ $u \tilde g(u) > 0$ for $u$ closed to $0$ and $u \neq 0$.\\ The converse is also true. When $0$ is an isochronous center of  $(E_c)$ this implies that $0$ is an isochronous center of $(E)$.\\
In fact, since \ $f$ and $g$ \ are analytic  \ $x g(x) > 0$\ in \ $N_0 - \{0\}$\ is a necessary and sufficient condition for the origin $0$ to be a center.\\   
So, We may assert that \ $0$\ is a center of \ $(E)$\ if and only if it is a center of \ $(E_c)$.\\

 In this case the integral \ $\tilde G(u)$\ may be expressed in terms of $x$
 $$\tilde G(u) = \int_0^x g(\sigma) e^{2F(\sigma)} d\sigma $$
 since \ $ds = e^{F(\sigma)}d\sigma .$\\
Furthermore, \ $X$\ which is such that \ $X^2 = 2\int_0^x g(s) e^{2F(s)} ds$ \ must verify  $$X \frac {dX}{du} = \tilde g(u).$$
On the other hand, by Lemma 2-2 (Urabe Theorem) \ $(E_c)$\ has an isochronous center at $0$ if and only if 
 $$\tilde g(u) = \frac {X}{1+h(X)} = g(x) e^{F(x)}$$
 where\  $h(X)$ \ is an odd analytic function such that $h(0) = 0$.\\
 Notice that by integration of \ $X \frac {dX}{du} = \tilde g(u) = \frac {x}{1 + h(X)}$\ one gets 
 $$u = X + H(X) = \int_0^x e^{F(s)} ds$$ 
 where \ $H(X) = \int_0^x h(\sigma) d\sigma$.\\
 In the isochronicity case \ $H(X)$\ must be even and must verify $H(0) = 0$.\\
 
Finally, the last part of Theorem 2-6  is a trivial consequence of the preceding. Indeed,  \ $h(X) \equiv 0$\ (or \ $H(X) \equiv 0$\ since $H$ is even)  means the equivalence reduces to  $$X = g(x) e^{F(x)} =  \phi (x) \qquad \Longleftrightarrow \qquad (E)\ {\mbox has \ an \ isochronous \ center}.$$
In the last case, $f$ and $g$ must be odd.

  \bigskip
{\bf Remark}\quad Theorem 2-6 characterizes all the isochronous centers of $(E)$ where $f$ and $g$ are analytic functions. When $f$ and $g$ are $C^1$ we get only necessary conditions. This characterization is based on the existence of an implicit function $h$ which is reduced to $h \equiv 0$ when $f$ and $g$ are odd. The calculation of the successive derivatives of\ $\tilde g(u) = \frac {X}{1 + h(X)}$\ allows us to establish an algorithm in the same manner of the one obtained thanks to the derivatives of the period function. This algorithm will permit us to find conditions for a critical point to be an isochronous center of (E).\\ 
  
  \bigskip
{\bf Corollary 2-7}\quad {\it Let $f,g$ be analytic function, then (E) has an isochronous center at  $0$ if and only if 
$$ \frac {dg}{dx} + f(x) g(x) = \frac {1+h(X) - h'(X) X}{(1+h(X))^3}$$ where $h(X)$ is an odd Urabe function.\\
In particular, when $f$ and $g$ are odd one gets the equivalence (in the isochronous case) $$\frac {dg}{dx} + f(x) g(x) =  1 \qquad \Longleftrightarrow \qquad h(X) \equiv 0$$}

  \bigskip
  
  Indeed, hypothesis \ $ \frac {dg}{dx} + f(x) g(x) = \frac {1+h(X) - h'(X) X}{(1+h(X))^3}$\ implies by integration \ $\frac {X}{1+h(X)} = g(x) e^{F(x)}$\ since \ $g(0) = 0$\ and \ $h(0) = 0.$\ Then by Theorem B, \ $0$\ is an isochronous center.\\
  The converse is trivial. \\
   More precisely, solving the differential equation $$ \frac {1+h(X) - h'(X) X}{(1+h(X))^3} = 1$$  one gets the solutions \ $1 + h(X) = \sqrt {\frac {X^2}{X^2+K}}$. For the boundary condition \ $h(0) = 0$ one have a unique solution :\ $h(X) \equiv 0$. \\
  Thus, Corollary 2-7 improves Lemma 2 and Corollary 9 of [S] which are limited to the case \ $ \frac {dg}{dx} + f(x) g(x) = 1$. \\
  
  Another consequence of Theorem 2-6 is the following which yields another characterization for isochronous centers of (E) and may be deduced from [C-J].
  
  \bigskip
     
    {\bf Corollary 2-8}\quad {\it Under the hypotheses of Theorem B if in addition :\\
    \ $ \frac {dg}{dx} + f(x) g(x) - 1$ \ is a polynomial in \ $\phi (x) = \int_0^xe^{-F(t)}dt $
 i.e. $$g(x) = e^{-F(x)} [\phi (x) + c_2(\phi (x))^2 + c_3(\phi (x))^3 +...+c_n(\phi (x))^n ]$$
then (E) has an isochronous center at $0$ if and only if \  $g(x) = \phi (x) e^{-F(x)}.$}

   \bigskip
   Indeed, hypothesis \ $g(x) = e^{-F(x)} [\phi (x) + c_2(\phi (x))^2 + c_3(\phi (x))^3 +...+c_n(\phi (x))^n ]$\ means 
 $$\tilde g(u) = g(x) e^{F(x)} = u + c_2 u^2 + c_3 u^3 +.....+c_n u^n.$$
 Then, following [C-J] Equation $(E_c)$ has isochronous center at $0$ and \ $\tilde g(u)$\ is a polynomial, then necessarily \ $\tilde g(u) = u $\ or equivalently \ $ g(x) = \phi (x) e^{-F(x)}.$
 
  \bigskip
  
  {\bf Corollary 2-9}\quad {\it Under the assumptions of Theorem 2-6 and suppose $0$ is an isochronous center of (E) then we have \ $X = g(x) e^{F(x)}= \int_0^x e^{F(s)} ds$\ if and only if \ $f(x)$\ and \ $g(x)$\ are odd functions.}
  
    \bigskip
 
  {\bf Remark} \quad Recall that when $0$ is an isochronous center of (E), hypothesis \ $f(x)$\ and \ $g(x)$\ odd implies \ $f(x) \equiv 0$\ and \ $g(x)$\ is linear, Corollary 10 of [S].\\
   When $0$ is an isochronous center of (E) and \ $h(X) \equiv 0$ then by Theorem 2-6,  \ $X = g(x) e^{F(x)},$\ i.e.   \ $X \equiv \phi (x).$\ Since  \ $u = X + h(X) \equiv X$\ it follows \ $g(x) e^{F(x)} = \phi (x)$.\ We get the converse if \ $f(x)$\ and \ $g(x)$\ are odd.\\

  \bigskip
  
  {\bf Corollary 2-10}\quad {\it Let $f,g$ be analytic functions and suppose \ $(E)$ \ has an isochronous center at the origin $0$. Let us define $$u = \int_0^x e^{F(s)} ds = X + H(X)$$ where $F(x) = \int_0^x f(s) ds,\ X$ is defined by $\frac {1}{2} X^2 = \int_0^x g(s) e^{2F(s)} ds,$ \\ $ \frac {dH}{dX}(X) = h(X)$\  and \ $H(0) = h(0) = 0$. Consider 
  $$\tilde g(u) = \frac {X}{1+\frac {dH}{dX}(X)} = g(x) e^{F(x)},$$   
 then $$\tilde g'(0) = 1,\quad \tilde g''(0) = g''(0) + f(0) = -3H''(0),$$ $$ \tilde g'''(0) = g'''(0) + 2f'(0) - f^2(0) = 15H''^2(0).$$}

 \bigskip

Recall here \ $\tilde g = \tilde g(u),\quad g = g(x)$\ and \ $H = H(X).$
 
{\bf Proof}\quad We calculate the successive derivatives of \ $\tilde g(u).$\ One gets:
$${(1+h)^{4}} \tilde g'(u) =   (-X-Xh) h'' +2\,X  (h'   
 ) ^{2}-2\,h'   -2\,  (h'h    )  $$ $$(1+h)^{6} \tilde g''(u) = -X(1+h)^{2} h'''+7\,  (-3/7\,h   -3/7+Xh'    )   (1+h) h''  -8\,
  (h')^{2}  (-1-h
   +Xh')$$
 $$\tilde g'''(u) =  [-X  (1+h    ) ^{3}{\frac {d^{4}}{d{X
}^{4}}}h   +11\,  (-4/11\,h   +X
h'-4/11 )   (1+h)^{2}h''' $$ $$  +7\,
X  (1+h    ) ^{2}  (h'')^{2}-59\,  (-{\frac {35}{59}}\,
h   +Xh'   -{\frac {35}{59
}} )   (1+h    )   h' h''$$ 
$$ +48\,  (h') ^{3}
  (-1-h   +Xh'h   
 )  ]   (1+h)^{-7}.$$
 Moreover, 
 $$\tilde g'(u) = \left( {\frac {d}{dx}}F \left( x \right)  \right) g \left( x \right) 
+{\frac {d}{dx}}g \left( x \right) ;$$
$$\tilde g''(u) =\left(  \left( {\frac {d^{2}}{d{x}^{2}}}F \left( x \right)  \right) g
 \left( x \right) + \left( {\frac {d}{dx}}F \left( x \right)  \right) 
{\frac {d}{dx}}g \left( x \right) +{\frac {d^{2}}{d{x}^{2}}}g \left( x
 \right)  \right) {e^{-F \left( x \right) }} ;$$

 $$ - \left( {e^{F(x)}} \right) ^{2}\tilde g'''(u) =  - \left( {
\frac {d^{3}}{d{x}^{3}}}F \left( x \right)  \right) g \left( x
 \right) -2\, \left( {\frac {d^{2}}{d{x}^{2}}}F \left( x \right) 
 \right) {\frac {d}{dx}}g \left( x \right)$$ $$-{\frac {d^{3}}{d{x}^{3}}}g
 \left( x \right) + \left( {\frac {d}{dx}}F \left( x \right)  \right) 
 \left( {\frac {d^{2}}{d{x}^{2}}}F \left( x \right)  \right) g \left( 
x \right) + \left( {\frac {d}{dx}}F \left( x \right)  \right) ^{2}{
\frac {d}{dx}}g \left( x \right)   .$$
We then deduce \ $\tilde g''(0) = g''(0) + f(0)$\ and \ $\tilde g'''(0) = g'''(0) + 2f'(0) - f^2(0).$

\bigskip

{\bf Remark}\quad By the same way one obtains
$$\tilde g^{(4)}(0) = g^{(4)}(0) - 2f(0) g'''(0) - f^2(0) g''(0) + 2f'(0)g''(0) - 7f(0)f'(0) + f''(0) + 2f^{3}(0) $$
$$\tilde g^{(5)}(0) = g^{(5)}(0) - 6f^{4}(0) + 5f^{3}(0)g''(0) + [29f'(0) + 5g'''(0)]f^{2}(0) + [-5g^{(4)}(0) - $$ 
$$ 19f''(0) 15f'(0)g''(0)]f(0) - 8f'^2(0) + 5f''(0)g''(0) + 4f'''(0) .$$
The proofs are omitted (the interested reader may used {\it Maple} for example).\\

Moreover, using again the relation between $\tilde g(u)$ and $H(X)$ one finds in the isochronous case
$$\tilde g^{(4)}(0) = 3780\,h'''(0) (h'(0))+945\,{(h'(0))}^{4},$$
$$\tilde g^{(5)}(0) = -1800\,h'''(0)-105\,{(h'(0))}^{3}.$$
This procedure of calculate different derivatives of $\tilde g $ and $H$ at $0$ allows us to define the algorithm described above for finding conditions for a critical point to be an isochronous center of $(E)$.

 \bigskip
  
  {\bf Corollary 2-11}\quad {\it Let $f,g$ be analytic functions and consider Equation $$(E) \qquad \ddot x + f(x)\dot x^2+g(x) = 0$$ with a center at the origin $0$.\ Let
  $$S(f,g) =  5g''^2(0) + 10g''(0)f(0) + 8f^2(0) - 3g'''(0) - 6f'(0).$$ 
 Then the following holds:\\
 (a)- $S(f,g) > 0$\ then the period function \ $T$\ increases in a neighborhood of $0$.\\ 
 (b)- $S(f,g) < 0$\ then the period function \ $T$\ decreases in a neighborhood of $0$.\\ 
 (c)- If (E) has an isochronous center at $0$ then \  $S(f,g) = 0.$}
 
 \bigskip
  
  Corollary 2-11 may be deduced from the preceding one. We may also use the Schaaf criteria  for the monotonicity of the period function for a conservative system, [Sc]. Indeed, 
  $$5 \tilde g''^2(0) - 3\tilde g'(0)\tilde g'''(0) = 5[g''(0) + f(0)]^2 - 3[g'''(0) + 2f'(0) - f^2(0)] $$ $$ \qquad = 5g''^2(0) + 10g''(0)f(0) + 8f^2(0) - 3g'''(0) - 6f'(0) = S(f,g).$$
  Notice that expression of \ $S(f,g)$\ (which corresponds to the coefficient of the first nonlinear term of the period function) may also be obtained by Proposition 2-4 since we have seen
  $$S(f,g) = p_2 = \frac {\pi}{2}\frac {d^3H}{dX^3}(0).$$
  When $f$ and $g$ are odd, $S(f,g)$ reduces to \ $S(f,g) = 8f^2(0) - 3g'''(0) - 6f'(0),$
  (see Corollary 6 of [S]).
  
  \newpage

\section{Applications to Loud systems}
In this part, we will apply preceding results to quadratic systems. In particular, the algorithm  presented above allows us to give another characterization of the isochrones of Loud systems in giving a new way of deriving the necessary condition. We then obtain a simple proof of Loud result for these systems.\\

After a rotation of coordinates, the Bautin's system may be transformed to  the form of the general  Loud system : 
$$(L_{B,D,F}) \cases{
   \dot x = - y + Bxy & \cr
   \dot y =   x + Dx^2 + Fy^2& \cr }
$$
for some $B,D$ and $F$ real parameters.\\
Notice that if the parameter $B \neq 0$, then another change of variables $u = Bx$ and $v = By$ the Loud system brings to the dehomogenized form
$$(L_{D,F}) \cases{
   \dot x = - y + xy & \cr
   \dot y =   x + Dx^2 + Fy^2.& \cr }
$$
Loud showed by direct integration of the systems that for the choice of the four pairs \ $(D,F)$
$$ I_1(0,1);\quad I_2(-\frac {1}{2},2);\quad I_3(0,(\frac {1}{4});\quad I_4(-\frac {1}{2},\frac {1}{2})$$
the corresponding system $(L_{D,F})$ has an isochronous center at $0$. Using Urabe theorem, Loud showed these are the only isochrones. \\
Recall that [C-J] proved the period coefficients \ $p_k,\ k\geq 2$\ for the dehomogenized Loud system are in the ideal \ $(p_2,p_4)$\ in the local ring \ $R(D,F)$\ localized at any of the isochrones \ $I_1, I_2, I_3, I_4$.\ Moreover, \ $p_2$\ is independant with respect to \ $p_4$\ at each isochrone.\\

Define a new independent time variable \ $\tau $\ by setting \ $t \rightarrow \tau$\ such that \ $dt = \psi(x) d\tau$\ where \ $\psi(x)$\ is \ $C^1$\ and \ $\psi(0) = 1$.

\bigskip

{\bf Lemma 3-1}\quad {\it By this change of time scale the dehomogenized Loud system \ $(L_{D,F})$\ is equivalent to a Lienard type equation (E)
$$(E_t)\qquad \psi (x) \dot x = y, \quad \psi (x) \dot y = - g(x) - f(x) y^2 $$
where \ $ f(x) = \frac {F+1}{1-x}-\frac {1}{\psi}\frac {d\psi }{dx}$\ and \ $g(x) = x(1-x)(1+Dx)\psi ^2 ,$\\
In particular, when \ $\psi \equiv 1$\ it is equivalent to $$(E_L)\qquad \ddot x +  \frac {F+1}{1-x} \dot x^2 + x(1 - x)(1 + Dx) = 0$$ } 

 \bigskip

 We shall prove this in two steps. Let  
 $$\frac {dx}{d\tau } = \dot x \psi, \quad  \frac {d^2x}{d\tau ^2} = \ddot x \psi ^2 + \dot x^2 \psi \frac {d\psi }{dx}$$ 
 $(L_{D,F})$ implies  $$\ddot x = - \dot y + \dot xy + x\dot y = \dot y(x - 1) + \dot xy$$ 
$$\ddot x =  (x + Dx^2 + Fy^2)(x - 1) + y^2(x - 1)$$ 
We then obtain 
$$\frac {d^2x}{d\tau ^2} = [(x + Dx^2 + Fy^2)(x - 1) + y^2(x - 1)]\psi ^2 + \dot x^2 \psi \frac {d\psi }{dx}$$
We thus deduce 
$$\frac {d^2x}{d\tau ^2} = [(x + Dx^2)(x - 1) + F(\frac {\dot x}{x-1})^2(x - 1) + (\frac {\dot x}{x-1})^2(x - 1)]\psi ^2 + \dot x^2 \psi \frac {d\psi }{dx}$$
which implies
$$\frac {d^2x}{d\tau ^2} = (x + Dx^2)(x - 1)\psi ^2 + [\frac {(F+1)\psi ^2}{x-1} + \psi \frac {\psi }{dx}] \dot x^2$$
$$ = (x + Dx^2)(x - 1)\psi ^2 + [\frac {(F+1)}{x-1} + \frac {1}{\psi } \frac {\psi }{dx}](\frac {dx}{d\tau})^2.$$

Notice that in changing the time scale : $t \rightarrow \tau , \ (L_{D,F})$\ is equivalent to a Lienard type system
$$(E_\tau ) \qquad \frac {dx}{d\tau} = y, \qquad \frac {dy}{d\tau} = - g(x) - f(x) y^2.$$ 
So, the systems $(E_t)$\ and \ $(E_\tau )$\ have the same phase portraits, but different period functions.\\
For the trivial case \ $\psi \equiv 0$,\ one gets the following Lienard type equation equivalent to  $(L_{D,F})$
$$(E_L)\qquad \ddot x +  \frac {F+1}{1-x} \dot x^2 + x(1 - x)(1 + Dx) = 0.$$ 
As a consequence of Theorem 2-6 and Corollary 2-10 one deduce the next result which yields a new way for obtained  necessary and sufficient condition for a dehomogenized Loud system $L_{D,F}$ to have an isochronous center at $0$. Our method differs of the one used by Loud himself [L] and Chicone and Jacobs [C-J]. Besides their method based on the vanishing of the {\it period quantity }\footnote{see the definition in the introduction} gives (in addition to the four isochrones \ $I_1, I_2, I_3, I_4$)\ the three weak centers denoted \ $L_1,L_2 \ and\ L_3.$ \ In fact, \ $p_2$\ and \ $p_4$\ have at most eight common zeros counted up to multiplicity. Opposite to our approach which yields exactly the four isochrones and no more. 

 \bigskip
 
{\bf Theorem 3-2}\quad {\it Let $X$ defined by $\frac {1}{2} X^2 = \int_0^x s(1+Ds)(1-s)^{-2F-1} ds$ \ and the function \ $H(X)$ is such that  $$  \frac {1}{F} [(1-x)^{-F} - 1]= X + H(X).$$
The dehomogenized Loud system $(L_{D,F})$ which is equivalent to the Lienard type equation
$$(E_L)\qquad \ddot x +  \frac {F+1}{1-x} \dot x^2 + x(1 - x)(1 + Dx) = 0$$ 
 has an isochronous center at the origin $0$ if and only if \ $H(X)$\ is even and
$$x (1-x)^{-F} (1+Dx) = \frac {X}{1+\frac {dH}{dX}}.$$
Moreover, \ $(L_{D,F})$\ has an isochronous center at $0$ for only the following values of the pairs}  \ $(D,F)$  :  
$$ I_1(0,1);\quad I_2(-\frac {1}{2},2);\quad I_3(0,\frac {1}{4});\quad I_4(-\frac {1}{2},\frac {1}{2})$$  

\bigskip

{\bf Corollary 3-3}\quad {\it The Loud system $(L_{D,F})$ has an isochronous center at the origin $0$ if and only if
$$(F-2)D x^2 + (2D+F-1) x + 1 = \frac {1+h(X) - h'(X) X}{(1+h(X))^3}$$ where $h(X) = \frac {dH}{dX}$ is an odd Urabe function.}

\bigskip

{\bf Proof of Theorem 3-2}\quad Let a variable defined above \ $u$ \ such that  $$\tilde g(u) = e^{F(x)} g(x)$$ where 
$$g(x) = x(1-x)(1+Dx), \qquad f(x) = \frac {F+1}{1-x}$$ 
$$u = \phi (x) = \frac {1}{F} [(1-x)^{-F} - 1].$$ Hence $$ \tilde g(u) = e^{F(x)} g(x) = x(1+Dx)(1-x)^{-2F}$$ 
So, by Theorem 2-6\ $(L_{D,F})$\ has an isochronous center at $0$ if and only if 
$$\tilde g(u) = \frac {X}{1+\frac {dH}{dX}}.$$
where \ $H = H(X)$\ is an even function.

In deriving, one obtains 
$$\tilde g'(u) = D(F-2)x^2 + (F+2D-1)x + 1$$ 
since \ $\frac {du}{dx} = (1 - x)^{-F-1}.$ \ By the same way we calculate 
$$\tilde g''(u) = [2D(F-2)x + F+2D-1](1-x)^{F+1}$$ 
$$\tilde g^{(3)}(u) = [2D(F-2)(1-x)^{F+1} - 2D(F-2)(F+1)x(1-x)^{F} -$$ $$\qquad (F+1)F+2D-1](1-x)^{F}](1-x)^{F+1}$$
$$\quad = [-2D(F+2)(F-2)x -6D-F^2+1](1-x)^{2F+1}$$
$$\tilde g^{(4)}(u) = -2D(F^2-4)(1-x)^{3F+2} +$$
$$ (2F+1)[-2D(F+2)(F-2)x - 6D-F^2+1](1-x)^{3F+1}$$
$$ \quad = [2D(F+2)(F-2)(2F+2)x + (6D+F^2-1)(2F+1) - 2D(F+2)(F-2)]$$

Thus, $$\tilde g'(0) = 1, \quad \tilde g''(0) = F + 2D - 1\quad 
\tilde g^{(3)}(0) = - F^2 - 6D + 1$$
$$\tilde g^{(4)}(0) = 2F^3+12DF-2DF^2+F^2-2F+14D-1 = (F+1)[2F^2 - 2DF + 14D - F - 1]$$
$$\tilde g^{(5)}(0) = (E + 1) (-6 F^3  + F^2  + 10 D F^2  + 4 F - 40 D F + 1 - 30 D),...$$

On the other hand, starting from \ $\tilde g(u) = \frac {X}{1+h(X)}$\ one gets
$$\tilde g''(0) = -3h'(0),\quad
\tilde g^{(3)}(0) = - 4h''(0) + 15 h'^2(0)$$
$$\tilde g^{(4)}(0) = -105 h'^3(0) + 45 h'(0) h''(0) - 5 h^{(3)}(0)$$
 $$\tilde g^{(5)}(0) = -6 h^{(4)}(0) + 70 h''^2(0) - 25 h'(0)h^{(3)}(0) - 330 h''(0)h'^2(0) + 105 h'^4(0).$$
 We then have
 $$4 h''(0) = F^2 + 6D - 1 + \frac {15}{9}(F +2D -1)^2 = (\frac {2}{3})[4 F^2 + 10 DF + 10 D^2 - D - 5F +1]$$
We have already seen that a necessary condition to have an isochronous center is \ $h''(0) = 0$\ which is equivalent to the Loud isochronicity condition [L]
$$(C_1)\quad 4 F^2 + 10 DF + 10 D^2 - D - 5F +1 = 0.$$

Using the above calculus, we find an isochronicity condition (different of those given by Loud) which permits to obtain by another way the four pairs of isochrone Loud systems \ $(L_{D,F}).$\\
 Recall at first since \ $h''(0) = 0$ then by Corollary 2-11
 $$ h^{(3)}(0) = -3h'^3(0)$$
 $$\tilde g^{(3)}(0) = 15 h'^2(0)$$
 $$\tilde g^{(4)}(0) = -105 h'^3(0) - 5 h^{(3)}(0) = -90 h'^3(0) = \frac {10}{3} (F+2D-1)^3 $$
 By idendification we get
 $$\tilde g^{(4)}(0) = 2F^3+12DF-2DF^2+F^2-2F+14D-1 = $$ $$\frac {10}{3} (F+2D-1)^3 = -2(F^2+6D-1)(F+2D-1)$$
 After simplification we obtain the new isochronicity condition
 
$$(C_2)\quad  4F^3 + 24DF + 24D^2 + 2DF^2 - F^2 - 4F - 2D + 1 = 0$$

Combined with \ $(C_1)$\ we may assert 

\bigskip

{\bf Lemma 3-4}\quad {\it The two equations 
$$(C_1)\quad 4 F^2 + 10 DF + 10 D^2 - D - 5F +1 = 0$$
$$(C_2)\quad 4F^3 + 24DF + 24D^2 + 2DF^2 - F^2 - 4F - 2D + 1 = 0$$
have only the following common real solutions  
$${D = 0, F = 1}; {D = -1/2, F = 2}; {D = 0, F = 1/4}; {F = 1/2, D = -1/2}$$    }

\bigskip
{\bf Proof } \quad
We may use classical computational method.\\
Write their resultants respectively of $D$ and $F$ 
$$R_1(D) = 864\,{D}^{2}+22176\,{D}^{4}+7536\,{D}^{3}+25920\,{D}^{5}+9600\,{D}^{6}$$
$$R_2(F) = -17280\,{F}^{3}+192+9000\,{F}^{2}-2160\,F-6480\,{F}^{5}+15768\,{F}^{4}
+960\,{F}^{6}$$
Resolve now\ $R_1(D) = 0$\ and \ $R_2(F) = 0.$\
The first equation gives real solutions
$$D = 0 , \qquad D = \frac {-1}{2}.$$
The second equation gives 
$$F = 1, \quad F = 2, \quad F = 1/4, \quad F = 1/2$$
Thus, one obtains exactly the four points.\\
In fact, thanks to {\it Maple} in solving \ $(C_1)$\ and \ $(C_2)$\ one obtains five pairs of solutions. 
The four real pairs
$$\{D = 0, F = 1\}; \{D = -1/2, F = 2\}; \{D = 0, F = 1/4\}; \{F = 1/2, D = -1/2\}$$
and a complex solution.\\

Finally,  two cases may occur. First if\ $h'(0)= 0,$ \ then $$ 2D+F-1 = 0.$$
It implies 
$$\tilde g'(u) = D(F-2)x^2 + (F+2D-1)x + 1 \equiv 1.$$ 
Then,  $$D(F-2) = F + 2D - 1 = 0 $$ and necessarely \ $(D,F) = (0,1)$\ or \ $(D,F) = (-\frac {1}{2}, 2).$\ Moreover, it implies \ $ h(X) \equiv 0$\ and \ $H(X) \equiv 0.$\  Thus, one find again the Loud isochrone systems\ $(L_{0,1})$\  and \ $(L_{-\frac {1}{2}, 2})$.\\

For the non trivial case $h(X) \neq 0 $.\ 
Let us denote \ $h'(0) = a \neq 0$.\ Then the preceding calculus gives  $$h^{(3)}(0) = -3 a^3, \ h^{(5)}(0) = 45 a^5, \ h^{(7)}(0) = -1575 a^7,....$$
That means the functions must take the following form
$$H(X) = \frac {1}{a}\sqrt{1+a^2X^2}, \qquad h(X) = \frac {a X}{\sqrt{1+a^2X^2}}$$

Notice that for \ $D = 0, \quad F = \frac {1}{4}$\ one finds \ $a = \frac {1}{4}.$\ 
In this case the corresponding Urabe function is
$$ h(X) =  \frac {X}{\sqrt {X^2+16}}$$
 For \ $D = -\frac {1}{2}, \quad F = \frac {1}{2},$\ one finds \ $a = \frac {1}{2}$\ and the corresponding Urabe function is  
$$h(X) =  \frac {X}{\sqrt {X^2+4}}.$$ 
 By our method we prove there are no other center candidate to be isochronous.
 \newpage
 
\section{Others monotonicity and isochronicity cases}
\subsection{On a generalization}
One of the most general equation that reduces to Equation $(E)$ is
$$ (E_q) \qquad \left\{ \begin{array}{c} \dot x = - \alpha (x)y \\
 \dot y = \beta (x) + \xi (x) y^2 \end{array}\right. $$
 with \ $\alpha , \beta , \xi $\ analytic functions in \ $N_0$\ a neighborhood of $0$.\\
 Equation $(E_q)$ is equivalent to $(E)$ with
 $$f(x) = \frac {\xi (x) - \alpha '(x)}{\alpha (x) } ,\qquad g(x) = \alpha (x) \beta (x) .$$
 We prove the following without need to suppose $f$ and $g$ odd.
 
 \bigskip

{\bf Corollary 4-1}\quad {\it Let \ $\alpha , \beta , \xi $\ be analytic functions, with \ $\alpha (x) > 0$\ and \ $x\beta (x) > 0$\ in \ $N_0$ \ a neighborhood of $0$.\ Let  $$X^2 = \int_0^x \frac {\beta (x)}{\alpha (x)} e^{2\int_0^x \frac {\xi (t)}{\alpha (t)} dt} dx.$$
A necessary and sufficient condition for the origin to be an isochronous center of $(E_q)$ is 
 $$\alpha (x) \beta '(x) + \xi (x) \beta (x)  = \frac {1+h(X) - h'(X) X}{(1+h(X))^3}$$ where $h(X)$ is an odd Urabe function.\\
In particular, one gets the equivalence $$\alpha (x) \beta '(x) + \xi (x) \beta (x) =  1 \qquad \Longleftrightarrow \qquad h(X) \equiv 0$$} 

 \bigskip
 
 {\bf  Proof}\quad This result follows from Theorem 2-6.
\\ Indeed, since \ $\alpha (x) > 0$\ condition  \ $x\beta (x) > 0$\ is equivalent to \ $x g(x) > 0$\ in \ $N_0.$ \ Moreover, 
$$f(x) g(x) + g'(x) = \alpha (x) \beta '(x) + \xi (x) \beta (x)$$
One has by notations of the preceding section $$F(x) = \int_0^x f(x) dx = \int_0^x \frac {\xi (t)}{\alpha (t) }dt - Log\alpha (x)$$ 
$$e^{F(x)} = \frac {1}{\alpha (x)} e^{\int_0^x \frac {\xi (t)}{\alpha (t)} dt}$$
$$u = \phi (x) = \int_0^x e^{F(x)} dx = \int_0^x \frac {1}{\alpha (x)} e^{\int_0^x \frac {\xi (t)}{\alpha (t)} dt}dx$$
Then we define \ $\tilde g(u) =  g(x)e^{F(x)} = \beta (x) e^{\int_0^x \frac {\xi (t)}{\alpha (t)} dt}$.\\ 
Moreover, a necessary and sufficient condition of isochronicity is 
$$\tilde g(u) = \frac {X}{1 + h(X)}$$ where \ $h(X)$\ is an odd Urabe function.
So, the condition holds by Corollary 2-1 and from the derivative
$$\frac {d\tilde g}{du} = \alpha (x) \beta '(x) + \xi (x) \beta (x) = \frac {d}{dX}[\frac {X}{1 + h(X)}].$$

\subsection {Reduced Kukles systems}
In this paragraph we apply preceding results to determine the monotonicity of the period function of reduced Kukles systems with a center at the origin.\\  
These systems correspond to second order differential equations and can be written as cubic systems under the form
 
$$ (K) \qquad \left\{ \begin{array}{c} \dot x = - y \\
 \dot y = x + a_1x^2 + a_2xy + a_3y^2 + a_4x^3 + a_5x^2y + a_6xy^2 \end{array}\right. $$
depending on parameters \ $a_i, \ i = 1,2,..,6.$\\
There are only four classes of reduced systems with a center, two of them are reversible systems.  In the case when the system symmetric with respect to the $x$-axis, [R-S-T] proved that \ $a_2 = a_5 = 0$\ is a necessary and sufficient condition for the reduced system to have a center. That is
$$ (K_0) \qquad \left\{\begin{array}{c} \dot x = - y\\
 \dot y = x + a_1x^2 + a_3y^2 + a_4x^3 + a_6xy^2. \end{array}\right. $$
They also obtain a non elementary first integral and the bifurcation diagram. Although the problem of finding general conditions for the Kukles systems to have a center still yet unsolved.\\

$(K_0)$ is also related to our equation 
$$(E) \qquad \ddot x + f(x)\dot x^2+g(x) = 0$$ with 
$$f(x) = a_3 + a_6x, \qquad g(x) = x + a_1x^2 + a_4x^3.$$
[S] interested in the monotonicity of the period function for system $(K_0)$ and found some sufficient condition but with restrictive hypotheses.\\

\newpage

{\bf Corollary 4-2}\quad {\it Let the expression 
$$S_{K_0} = 10a_1^2+10a_1a_3+4a_3^2-9a_4-6a_6$$
(i) - If\ $S_{K_0} > 0$\ then the period function of $(K_0)$ is increasing at $0.$\\
(ii) - If\ $S_{K_0} < 0$\ then the period function of $(K_0)$ is decreasing at $0$.}

\bigskip
Indeed, this statement follows from Corollary 2-11 since
$$S(f,g) =  S_{K_0} = 5g''^2(0) + 10g''(0)f(0) + 8f^2(0) - 3g'''(0) - 6f'(0).$$ 
$$ S_{K_0} = 10a_1^2+10a_1a_3+4a_3^2-9a_4-6a_6.$$

{\bf Remark}\quad What it happens when \ $ S_{K_0} = 0$ ?
It is wellknown the origin $0$ is never isochronous (except for the linear case) since it is a weak center of order at most 3, [R-S-T]. Moreover, a such center has a perturbation with at most 3 local critical periods. \\
We can see that in another manner thanks to {\it Maple}. Let us take again the proof of Corollaries 2-10 and 2-11. We get from the relation which connects $\tilde g(u)$ and $H(X)$ 
$$\tilde g^{(4)}(0) = g^{(4)}(0) - 2f(0) g'''(0) - f^2(0) g''(0) + 2f'(0)g''(0) - $$
$$7f(0)f'(0) + f''(0) + 2f^{3}(0) = - 90H''^3(0)$$
and $$\tilde g''(0) = g''(0) + f(0) = -3H''(0)$$
Recall that \ $\tilde g(x) = g(x) e^{F(x)}$.\ After replacing we then obtain
$$\tilde g''(0) = 2a_1 + a_3 = -3H''(0)$$
$$\tilde g^{(4)}(0) = 24\,{\it a_4}\,{\it a_3}+12\,{\it a_1}\,{\it a_6}+12\,{\it a_1}\,{{\it a_3}}
^{2}+15\,{\it a_6}\,{\it a_3}+3\,{{\it a_3}}^{3}-$$
$$3\, \left( 6\,{\it a_4}+6
\,{\it a_1}\,{\it a_3}+3\,{\it a_6}+3\,{{\it a_3}}^{2} \right) {\it a_3}-3
\, \left( 2\,{\it a_1}+2\,{\it a_3} \right) {\it a_6}+3\, \left( 2\,{\it 
a_1}+2\,{\it a_3} \right) {{\it a_3}}^{2}-$$
$$2\, \left( 6\,{\it a_4}+6\,{\it 
a_1}\,{\it a_3}+2\,{\it a_6}+4\,{{\it a_3}}^{2}-2\, \left( 2\,{\it a_1}+2\,
{\it a_3} \right) {\it a_3} \right) {\it a_3}- \left( 2\,{\it a_1}+{\it a_3
} \right) {\it a_6}+ \left( 2\,{\it a_1}+{\it a_3} \right) {{\it a_3}}^{2}
- $$ 
$$\left( 6\,{\it a_4}+6\,{\it a_1}\,{\it a_3}+2\,{\it a_6}+4\,{{\it a_3}}^{
2}-2\, \left( 2\,{\it a_1}+2\,{\it a_3} \right) {\it a_3}- \left( 2\,{
\it a_1}+{\it a_3} \right) {\it a_3} \right) {\it a_3}-10/3\, \left( 2\,{
\it a_1}+{\it a_3} \right) ^{3}$$
$$= - 90H''^3(0) = \frac {10}{3}(2a_1 + a_3)^3$$
Simplifying the above expression one gets
$$\Sigma_{K_02} = -\frac {4}{3}\,{{\it a_3}}^{3}-22\,{\it a_1}\,{{\it a_3}}^{2}+\frac {1}{3}\, \left( -120\,{
{\it a_1}}^{2}-36\,{\it a_4}-21\,{\it a_6} \right) {\it a_3}+4\,{\it a_1}\,
{\it a_6}-{\frac {80}{3}}\,{{\it a_1}}^{3} = 0 $$
By the same way in considering \ $\tilde g^{(5)}(0)$\ and its connection with derivatives of the function \ $H(X)$\ 
$$\tilde g^{(5)}(0) - 630H''^4(0) = -4\,{{\it a_3}}^{4}+8\,{\it a_1}\,{{\it a_3}}^{3}+ \left( 22\,{\it a_6}+18
\,{\it a_4} \right) {{\it a_3}}^{2}-26\,{\it a_1}\,{\it a_6}\,{\it a_3}-$$ 
$$\qquad 8\,{{\it a_6}}^{2} - \frac {70}{9}(2a_1 + a_3)^4$$
One obtains another relation
$$\Sigma_{K_03} = -4\,{{\it a_3}}^{4}+\frac {1}{9}\, \left( 72\,{\it a_1}-70 \right) {{\it a_3}}^{3}
+1/9\, \left( -420\,{\it a_1}+198\,{\it a_6}+162\,{\it a_4} \right) {{
\it a_3}}^{2}+$$
$$\frac {1}{9}\, \left( -234\,{\it a_1}\,{\it a_6}-840\,{{\it a_1}}^{2}
 \right) {\it a_3}-8\,{{\it a_6}}^{2}-{\frac {560}{9}}\,{{\it a_1}}^{3}.$$
 
 Solving now \ $S_{K_0} = 0,\quad \Sigma_{K_02} = 0$,\ one gets the following \\
 (i) if \ $a_1 a_3 \neq 0$ and $-4a_1+3a_3 \neq 0$ then 
 
 $$a_6 = -(2/3)\frac {(53a_1a_3^2+40a_1^3+10a_3^3+80a_1^2a_3)}{(-4a_1+3a_3)},$$
 $$ a_4 = (2/9)\frac {(20a_1^3+75a_1^2a_3+60*a_1a_3^2+16a_3^3)}{(-4a_1+3a_3)}$$
 (ii) if \ $a_1 = 0, a_3 = 0$ then
 $$ a_4 = -\frac {1}{3}a_6.$$
Moreover, the system of three equations 
$$S_{K_0} = 0,\quad \Sigma_{K_02} = 0,\quad \Sigma_{K_03} = 0$$
has only the trivial solution \ $a_1 = a_3 = a_4 = a_6 = 0$\  corresponding to the linear isochrone.

\subsection{A cubic system}
Cubic systems with non homogeneous singularities take the following form
$$\left\{ \begin{array}{c} \dot x = -y + b_1x^2 + b_2y^2 + b_3xy + b_4x^3 + b_5x^2y + b_6xy^2 + b_7y^3\\
\dot y = x + a_1x^2 + a_2y^2 + a_3xy + a_4x^3 + a_5 yx^2 + a_6xy^2 + a_7y^3
\end{array}\right. $$
In polar coordinates this system may be written under the form

$$(C*) \qquad \left\{ \begin{array}{c} \dot r = r^2(R_3 \sin 3\theta + R_1 \sin \theta )+ r_3(R_4 \sin 4\theta + R_2 \sin 2\theta )\\
\dot \theta = 1 + r(R_3 \cos 3\theta + r_1 \cos \theta ) + r^2(R_4 \cos 4\theta + R_2 \cos 2\theta + r_0) \end{array}\right. $$
Pleshkan studied systems with homogeneous singularities and an isochronous center at $0$. He proved that there are only four different classes of such systems.\\
Chavarriga and Garcia ([C-S] sect. 12) considered cubic reversible systems of the form above $(C*)$.\\
They completely classified the case $R_3 = 0.$ More precisely, they proved that such systems with an isochronous center at $0$ and verifying $R_3 = 0$ and $ R_4 \neq 0$ can be brought to one of the Pleshkan cubic homogeneous systems, denoted by $(S_1^*), (S_2^*), (S_3^*) \ or \ (S_4^*)$ (in using their terminology). The remainig case \ $R_3 \neq 0$ still open. There are only few examples of such systems verifying $R_3 \neq 0$. \\
Garcia  considered families of isochrone reversible cubic systems of the form (see [G] p 108)

$$   \left\{ \begin{array}{c} \dot x = - y + (2-r_1+R_1)xy + (3 R_4+2 R_2-r_2-r_0)x^2y + (r_2-r_0-R_4)y^3\\
 \dot y = x + (1+r_1)x^2 + (R_1-1)y^2 + (R_4+r_0+r_2)x^3 + (r_0+r_2+2 R_2-3 R_4)x y^2 \end{array}\right. $$
where $R_1, R_2, R_3, R_4, r_0, r_1$ and $r_2$ are the coefficients defined in $(C*)$ satisfying the condition \ $r_2 - r_0 - R_4 = 0.$

Let us consider the following cubic system depending on five parameters. We will give necessary and sufficients conditions so that $(C)$ has an isochronous center 
$$ (C) \qquad \left\{ \begin{array}{c} \dot x = - y + b x^2y\\
 \dot y = x + a_1x^2 + a_3y^2 + a_4x^3 + a_6xy^2 \end{array}\right. $$
depending on real parameters \ $b, a_1, a_3, a_4, a_6.$ \\
Here  $$4 R_3 = a_1 - a_3.$$ 
This system having a center at the origin $0$. Moreover, we establish necessary and sufficient conditions so that this center is isochronous and no additional condition on these coefficients is required.  \\

More precisely, using our algorithm above we shall prove that System $(C)$ possesses only four classes of isochrone systems. In particular, we produce examples of cubic reversible systems verifying $R_3 \neq 0$ with an isochronous center at $0$.
The following result improves the study of cubic reversible systems started by Garcia, [G], p 108. 

\bigskip

{\bf Theorem 4-3}\quad {\it Let us consider the expression 
$$S_C = 20a_1^2+20a_1a_3+8a_3^2-18a_4-6a_6+6b$$
(i) - If\ $S_{C} > 0$\ then the period function of $(C)$ is increasing at $0.$\\
(ii) - If\ $S_{C} < 0$\ then the period function of $(C)$ is decreasing at $0$.\\
(iii) - The system $(C)$ has an isochronous center at $0$ if and only if it reduces to the one of the following}
$$(I)\quad {a_4 = -(2/3)b, a_1 = 0, a_3 = 0, a_6 = 3b, b = b}$$
 $$(II)\quad  {a_1 = 0, a_3 = 0, a_6 = b, a_4 = 0, b = b}$$
 $$(III)\quad  {a_4 = (1/14)a_3^2, a_6 = (3/7)a_3^2, b = (1/7)a_3^2, a_1 = -(1/2)a_3, a_3 = a_3}$$
 $$(IV)\quad  {a_6 = a_3^2, a_4 = 0, b = (1/2)a_3^2, a_1 = -(1/2)a_3, a_3 = a_3}.$$

\bigskip

{\bf Proof}\quad The system $(C)$ is equivalent to Equation 
$$(E) \qquad \ddot x + f(x)\dot x^2+g(x) = 0$$ with 
$$f(x) = \frac{a_3 + a_6x + 2bx}{1 - bx^2}, \qquad g(x) = (x + a_1x^2 + a_4x^3)(1 - bx^2).$$
We will apply preceding methods of Corollaries 4 and 5.\\
Indeed, a calculus gives the following
$$ e^{F(x)} = {e^{-1/2\,{\frac {\ln  \left( -1+b{x}^{2} \right) {\it a_6}}{b}}-\ln 
 \left( -1+b{x}^{2} \right) +{\frac {{\it a_3}\,{\it arctanh} \left( 
\sqrt {b}x \right) }{\sqrt {b}}}}}$$
$$\tilde g(u) =  \left( x+{\it a_1}\,{x}^{2}+{\it a_4}\,{x}^{3} \right)  \left( 1-b{x}^{
2} \right) {e^{-1/2\,{\frac {\ln  \left( -1+b{x}^{2} \right) {\it a_6}}
{b}}-\ln  \left( -1+b{x}^{2} \right) +{\frac {{\it a_3}\,{\it arctanh}
 \left( \sqrt {b}x \right) }{\sqrt {b}}}}}$$
 where \ $\frac {du}{dx} = e^{F(x)}.$
 $$\tilde g(u) = - \left( x+{\it a_1}\,{x}^{2}+{\it a_4}\,{x}^{3} \right){e^{-1/2\,{\frac {\ln  \left( -1+b{x}^{2} \right) {\it a_6}}
{b} +{\frac {{\it a_3}\,{\it arctanh}
 \left( \sqrt {b}x \right) }{\sqrt {b}}}}}}$$
 Recall by Corollary 2-11
 $$S_C = 5 \tilde g''^2(0) - 3\tilde g'(0)\tilde g'''(0) = 5[g''(0) + f(0)]^2 - 3[g'''(0) + 2f'(0) - f^2(0)] $$ $$ \qquad = 5g''^2(0) + 10g''(0)f(0) + 8f^2(0) - 3g'''(0) - 6f'(0)$$

So, we evaluate the expression 
$$S_C = A =\left( 20\,{{\it a_1}}^{2}+20\,{\it a_1}\,{\it a_3}+8\,{{\it a_3}
}^{2}-18\,{\it a_4}-6\,{\it a_6}+6\,b \right)$$

 We thus prove part (i) and (ii) of Theorem 4-3.\\
 Now, in order to prove part (iii), we need to calculate \ the next derivatives of $\tilde g(u).$\\ It yields respectively after simplication and thanks to {\it Maple}
 $$\tilde g^{(4)}(0) = -16\,b{\it a1}-12\,{\it a3}\,{\it a4}+4\,{\it a3}\,b+4\,{\it a1}\,{
\it a6}-2\,{\it a1}\,{{\it a3}}^{2}-7\,{\it a3}\,{\it a6}+2\,{{\it a3}
}^{3}$$
 
 $$\tilde g^{(5)}(0) = -120\,b{\it a4}+80\,b{\it a1}\,{\it a3}+16\,{b}^{2}-8\,b{\it a6}+30\,{
{\it a3}}^{2}{\it a4}$$ $$-8\,{{\it a6}}^{2}-10\,{{\it a3}}^{2}b+10\,{{\it 
a3}}^{3}{\it a1}+29\,{{\it a3}}^{2}{\it a6}-6\,{{\it a3}}^{4}-30\,{
\it a3}\,{\it a1}\,{\it a6},$$
and so on.\\  
On the other hand, another calculation yields (see also the remark following Corollary 2-10)                                        
$$\tilde g^{(4)}(0) = 3780\,h'''(0) (h'(0))+945\,{(h'(0))}^{4},$$
$$\tilde g^{(5)}(0) = -1800\,h'''(0)-105\,{(h'(0))}^{3},...$$
After identification and thanks to Maple we find only four (real) solutions (the details will be given in Appendice 2)  
 $$(I)\quad {a_4 = -(2/3)b, a_1 = 0, a_3 = 0, a_6 = 3b, b = b}$$
 $$(II)\quad  {a_1 = 0, a_3 = 0, a_6 = b, a_4 = 0, b = b}$$
 $$(III)\quad  {a_4 = (1/14)a_3^2, a_6 = (3/7)a_3^2, b = (1/7)a_3^2, a_1 = -(1/2)a_3, a_3 = a_3}$$
 $$(IV)\quad  {a_6 = a_3^2, a_4 = 0, b = (1/2)a_3^2, a_1 = -(1/2)a_3, a_3 = a_3}$$
 
 {\bf Case} $(I)$. Indeed, the following system 
 $$ (C_I) \qquad \left\{ \begin{array}{c} \dot x = - y + \frac {1}{3}a_6 x^2y \\ 
 \dot y = x - \frac {2}{9}a_6 x^3 + a_6 xy^2 \end{array}\right. $$ 
 corresponds to Equation (E) of the form
 $$\ddot x + \frac {5a_6 x}{3 - a_6 x^2} \dot x^2 + (x - \frac {2}{9}a_6 x^3)(1 - \frac {a_6}{3} x^2) = 0.$$
This equation has an isochronous center at $0$ since  the corresponding functions 
$$f(x) = \frac {5a_6 x}{3 - a_6 x^2}, \quad g(x) = (x - \frac {2}{9}a_6 x^3)(1 - \frac {a_6}{3} x^2)$$ are odd and are such that $$f(x)g(x) + g'(x) = \frac {5}{3}a_6 x (x - \frac {2}{9}a_6 x^3) + \left( 1-
\frac {2}{3}\,{\it a_6}\,{x}^{2} \right)  \left( 1-\frac {1}{3}\,{\it a_6}\,
{x}^{2} \right)$$
$$ -\frac {2}{3}\, \left( x-\frac {2}{9}\,{\it a_6}\,{x}^{3} \right) {\it 
a _6}\,x$$
$$=\frac {5}{3}a_6 x (x - \frac {2}{9}a_6 x^3) 1-\frac {5}{3}\,{\it a_6}\,{x}^{2}+{\frac {10}{27}}\,{{\it a_6}}^{2}{x}^{4} = 1$$
Then, by Corollary 2-7  System $(C_I)$ is isochrone. 
Notice that  $(C_I)$ is denoted by $(S_3^*)$ in the classification of isochronous homogenous cubic systems of Pleshkan.\\

 {\bf Case} $(II)$. In this case $(C)$ reduces to 
 $$ (C_{II}) \qquad \left\{ \begin{array}{c} \dot x = - y + a_6 x^2y \\ 
 \dot y = x + a_6 xy^2 \end{array}\right. $$
 which corresponds to Equation (E) of the form
 $$\ddot x + \frac {3a_6 x}{1 - a_6 x^2} \dot x^2 + (x - a_6 x^3) = 0.$$
This equation has an isochronous center at $0$ since the corresponding functions  
$$f(x) = \frac {3a_6 x}{1 - a_6 x^2}, \quad g(x) = (x - a_6 x^3)$$ are odd and  $$f(x)g(x) + g'(x) = 3a_6 x^2 - 3a_6 x^2 + 1 = 1,$$
 then by Corollary 2-7  $(C_{II})$ is isochrone. \\
 In fact, $(C_{II})$ is a (trivial) cubic reversible degenerated system, i.e. $R_3 = R_4 = 0,$ (see Theorem 8.11 of [G]). Moreover, the change $z = y^2$ transforms $(C_{II})$ into the linear differential equation 
 $$(-1 + a_6 x^2) \frac {dz}{dx} = 2x + 2a_6 x z.$$
 It gives immediately a rational first integral and a transversal commuting system.\\
 
 The remaining cases are more instructive because they satisfy the condition $R_3 \neq 0$ ( $R_3 = 0$ being completely solved for cubic polynomial systems, [G]).\\  
 
 {\bf Case} $(III)$. In this case $(C)$ reduces to 
 $$ (C_{III}) \qquad \left\{ \begin{array}{c} \dot x = - y + \frac {1}{7}a_3^2 x^2y\\
 \dot y = x + \frac {-1}{2}a_3 x^2 + a_3y^2 + \frac {1}{14}a_3^2 x^3 + \frac {3}{7}a_3^2 xy^2 \end{array}\right. $$
$(C_{III})$ corresponds to Equation (E) of the form
$$\ddot x + \frac {a_3 + (5/7)a_3^2 x}{1 - (1/7)a_3^2 x} \dot x^2 + (x - (1/2)a_3 x^2 + (1/14)a_3^2 x^3){1 - (1/7)a_3^2 x} = 0.$$
Following Chavarriga this system has an isochronous center at $0$.\\

{\bf Case} $(IV)$. In this case $(C)$ reduces to 
 $$ (C_{IV}) \qquad \left\{ \begin{array}{c} \dot x = - y + \frac {1}{2}a_3^2 x^2y\\
 \dot y = x + \frac {-1}{2}a_3 x^2 + a_3y^2 + a_3^2 xy^2 \end{array}\right. $$
$(C_{IV})$ corresponds to Equation (E) of the form
$$\ddot x + \frac {a_3 + 2a_3^2 x}{1 - (1/2)a_3^2 x} \dot x^2 + (x - (1/2)a_3 x^2)(1 - (1/2)a_3^2 x) = 0.$$
Following Garcia (family A, p.109) $(C_{IV})$ posseses two invariant algebraic curves and a first integral and a linearized change of variables. Thus, $(C_{IV})$ is an isochrone system.

\bigskip

{\bf Corollary 4-4} \quad  {\it Under the hypothesis $R_3 \neq 0$ (i.e. $a_1 \neq a_3$) the system 
$$ (C) \qquad \left\{ \begin{array}{c} \dot x = - y + b x^2y\\
 \dot y = x + a_1x^2 + a_3y^2 + a_4x^3 + a_6xy^2 \end{array}\right. $$
 has an isochronous center at $0$ if and only it can be reduced to the one of the following} 
 $$ (C_{III}) \qquad \left\{ \begin{array}{c} \dot x = - y + \frac {1}{7}a_3^2 x^2y\\
 \dot y = x + \frac {-1}{2}a_3 x^2 + a_3y^2 + \frac {1}{14}a_3^2 x^3 + \frac {3}{7}a_3^2 xy^2 \end{array}\right. $$
 $$ (C_{IV}) \qquad \left\{ \begin{array}{c} \dot x = - y + \frac {1}{2}a_3^2 x^2y\\
 \dot y = x + \frac {-1}{2}a_3 x^2 + a_3y^2 + a_3^2 xy^2 \end{array}\right. $$

\bigskip

\newpage
\section{Appendice 1}
 
 Equation $(E)$
$$\qquad \ddot x + f(x)\dot x^2 + g(x) = 0$$  can be solved by use firstly a reduction of the order by the change  $$\dot x=p,\qquad \dot x=p\frac {dp}{dx}.$$
 The obtained first order equation can be solved in using the integrating factor $\mu (x)=e^{2\int_0^x f(s) ds}$.\\ Let us consider the case $$f(x) = \frac {- \lambda x}{1+\lambda x^2}, \qquad  g(x) = \frac {\alpha ^2 x}{1+\lambda x^2}$$  then \ $\mu(x) = \frac {1}{1+\lambda x^2}$.\\ Equation $(E)$ may be deduced from the associated Lagrangian which is $$L=\frac {1}{2}(\frac {1}{1+\lambda x^2})(x'^2-\alpha^2x^2).$$ The Lagrangian density \
${\cal L} = \frac {1}{2}(\frac {1}{1+\lambda \phi ^2})(\partial_\mu \phi \partial^\mu \phi - m^2\phi ^2)$\ 
 appears in some models of quantum field theory.\\
The one dimensional Schrodinger equation involving the potential $\frac {x^2}{1+\lambda x^2}$ was considered as an  example of anharmonic oscillator, see [L-R] for more details.\\
In order to know the growth of the period function we will apply Lemma 1 which transforms  $(E)$ 
$$\ddot x + \frac {- \lambda x}{1+\lambda x^2}\dot x^2 + \frac {\alpha ^2 x}{1+\lambda x^2} = 0$$
to a conservative system 
$$\dot u = v,\qquad \dot v = - \tilde g(u)$$
where $$\tilde g(u) = \alpha ^2 \frac {\sinh u}{(\cosh u)^3}.$$
On the other hand, since the functions $f$ and $g$ are odd (and as well as  $\tilde g (u)$) we may apply Corollary 2-5 to assert that 
$$ x[g(x)\phi'(x) \phi (x)g'(x) - \phi (x)g(x)f(x)]$$
has a minimum at $0$ if and only if $T$ is an increasing function.
\newpage

\section{Appendice 2: Complement of the proof of Theorem 4-3}
Let the functions $$f(x) = \frac{a_3 + a_6x + 2bx}{1 - bx^2}, \qquad g(x) = (x + a_1x^2 + a_4x^3)(1 - bx^2)$$
where \ $\tilde g(u) = \frac {X}{1 + h(X)}.$
The odd function $h(X)$ may be written
$$h(X) = d X + e X^3 + c X^5 + k X^7 +...$$ One find the two next derivatives of $\tilde g(u)$
$$ \tilde g^{(6)}(0) = -144\,{\it a_3}\,{b}^{2}+1080\,b{\it a_4}\,{\it a_3}+272\,{b}^{2}{\it a_1}
-46\,{\it a_1}\,{{\it a_6}}^{2}-90\,{{\it a_3}}^{3}{\it a_4}$$ $$-30\,{\it a_3}
\,{\it a_4}\,{\it a_6}+44\,b{\it a_1}\,{\it a_6}-400\,{{\it a_3}}^{2}b{\it 
a_1}+208\,{{\it a_3}}^{2}{\it a_1}\,{\it a_6}+92\,{\it a_3}\,b{\it a_6}$$ $$+24\,
{{\it a_3}}^{5}+97\,{\it a_3}\,{{\it a_6}}^{2}+30\,{{\it a_3}}^{3}b-52\,{{
\it a_3}}^{4}{\it a_1}-146\,{{\it a_3}}^{3}a $$
$$ \tilde g^{(7)}(0) = -1540\,{{\it a_3}}^{3}{\it a1}\,{\it a_6}+2240\,{{\it a_3}}^{3}b{\it a_1}+
308\,{{\it a_3}}^{5}{\it a_1}+294\,{{\it a_3}}^{4}{\it a_4}-98\,{{\it a_3}}
^{4}b$$ $$-272\,{b}^{3}-3808\,{b}^{2}{\it a_1}\,{\it a_3}-560\,{\it a_3}\,b{
\it a_1}\,{\it a_6}+874\,{{\it a_3}}^{4}{\it a_6}+504\,{{\it a_3}}^{2}{\it 
a_4}\,{\it a_6}$$ $$-8400\,b{\it a_4}\,{{\it a_3}}^{2}-880\,{{\it a_3}}^{2}b{
\it a_6}+840\,{\it a_3}\,{\it a_1}\,{{\it a_6}}^{2}-120\,{{\it a_3}}^{6}+
104\,{{\it a_6}}^{3}+1120\,{{\it a_3}}^{2}{b}^{2}$$ $$-969\,{{\it a_3}}^{2}{{
\it a_6}}^{2}-72\,{b}^{2}{\it a_6}+3696\,{b}^{2}{\it a_4}-168\,{\it a_4}\,
{{\it a_6}}^{2}+240\,b{{\it a_6}}^{2}+1512\,b{\it a_4}\,{\it a_6}$$
On the other hand, the relation between $\tilde g(u)$  and $h(X)$ gives
$$ \tilde g^{(6)}(0) = -840\,c-11340\,e{d}^{2}-10395\,{d}^{5}$$
$$ \tilde g^{(7)}(0) = 30240\,dc+135135\,{d}^{6}+207900\,{d}^{3}e+11340\,{e}^{2}.$$
Finally, the identification yields the four solutions 

$$(C_{I}) \  a_6 = 3 b, a_4 = -(2/3)b, a_1 = 0, a_3 = 0, \quad h(X) \equiv 0$$

$$(C_{II}) \ { a_6 = b, a_1 = 0, a_3 = 0, a_4 = 0, k = 0}, \quad h(X) \equiv 0$$

$$(C_{III})\ { b = (1/7)a_3^2, a_6 = (3/7)a_3^2, a_1 = -(1/2)a_3, a_4 = (1/14)a_3^2}, \quad 
h(X) = (1/3087)a_3^7 X^7 + ...$$

$$(C_{IV})\ {  a_1 = -(1/2)a_3, a_6 = a_3^2, b = (1/2)a_3^2, a_4 = 0}, \quad 
h(X) = (1/72)a_3^7 X^7 + ...$$

\newpage

{\bf REFERENCES}\\

\smallskip
[C-S] \ J. Chavarriga and M. Sabatini \quad {\it A survey of isochronous centers }\quad  Qual. Theory of Dyn. Systems {\bf vol 1} p.1-70, (1999).

\smallskip
{[C-J]}\ C. Chicone and M. Jacobs,
            {\it Bifurcation of critical periods for plane vector fields},\ 
             Trans. A.M.S. \ {\bf 312}(2), (1989), p.433-486.  
             
\smallskip
[G] I. Garcia, \quad {\it Contribution to the qualitative study of planar vector fields}\quad Dept. de Matematica. University of Lleida, (2000).            

\smallskip
[L-R] M. Lakshmanan and S. Rajasekar, \quad {\it Nonlinear dynamics. Integrability, Chaos and Patterns}, \ Advanced Texts in Physics, Springer-Verlag, Berlin, (2003).

\smallskip
[L] \ W.S. Loud \quad {\it The behavior of the period of solutions of certain plane autonomous systems near centers  }\quad Contr. Differential Equations, {\bf 3}, p. 21-36, (1964). 
             
\smallskip
[R-S-T] C.Rousseau, D.Schlomiuk, P.Thibaudeau, \quad {\it The centres in the reduced Kukles system}, \  Nonlinearity  {\bf 8}  no. 4, p.541--569, (1995).

\smallskip
 [S] M. Sabatini, \quad {\it On the period function of $x''+f(x)x'^2+g(x)=0$ }, \ J. Diff. Eq., 152, p. 1-18, 197 (2003).
 
 \smallskip
[Sc] \ R.Schaaf \quad {\it A class of Hamiltonian systems with increasing periods,} \quad J. Reine Angew. Math., {\bf 363},  (1985), 96-109.          

\smallskip
[U1]  M. Urabe \quad  {\it The potential force yielding a periodic motion whose period is an arbitrary continuous function of the amplitude of the velocity}\ Arch. Ration. Mech. Anal.,{\bf 11}, p.27-33, (1962).

\end{document}